\documentclass[11pt]{article}
\usepackage{amsfonts,latexsym,amsmath,amscd,geometry}
\usepackage{amssymb}
\usepackage{latexsym}

\newcommand \nc{\newcommand}
\newtheorem{theorem}{Theorem}[section]
\newtheorem{lemma}[theorem]{Lemma}
\newtheorem{proposition}[theorem]{Proposition}
\newtheorem{corollary}[theorem]{Corollary}
\newtheorem{definition}[theorem]{Definition}

\nc{\ba}{\begin{array}}\nc{\ea}{\end{array}}
\nc{\be}{\begin{eqnarray}}\nc{\ee}{\end{eqnarray}}
\nc{\beq}{\begin{equation}}\nc{\eeq}{\end{equation}}
\nc{\bex}{\begin{eqnarray*}}\nc{\eex}{\end{eqnarray*}}
\nc{\btm}{\begin{theorem}} \nc{\etm}{\end{theorem}}
\nc{\blm}{\begin{lemma}} \nc{\elm}{\end{lemma}}
\nc{\R}{\mathbb{R}} \nc{\va}{\varepsilon} \nc{\ls}{\limits}
\def\P{\Phi}\def\D{\Delta}\def\T{\mathbf{T}}\def\B{\mathbf{B}}
\def\pf{\noindent{\bf Proof.\quad}}\def\endpf{\hfill$\Box$}
\def\Li{L^{\infty}}\def\u{\hat{u}_0} \textwidth=150mm \textheight=225mm
\nc \qed {\hfill $\Box$}

\begin{document}
\title{Well-posedness for the heat flow of
polyharmonic maps with rough initial data}
\author{Tao Huang\footnote{Department of Mathematics, University of Kentucky,
Lexington, KY 40506} \quad Changyou Wang$^*$ }
\date{}
\maketitle

\begin{abstract} We establish both local and  global well-posedness of the heat
flow of polyharmonic maps from $\mathbb{R}^n$ to a compact
Riemannian manifold without boundary for initial data with small BMO norms.
\end{abstract}

\vskip 5mm
\section {Introduction}
\setcounter{equation}{0}
\setcounter{theorem}{0}

For $k\geq 1$, let $N$ be a $k$-dimensional compact Riemannian
manifold without boundary, isometrically embedded in some Euclidean
space $\mathbb{R}^l$. For $n\geq 2$ and $m\geq 1$, we consider the
$m$-th order energy functional
\bex
E_m(u)=\frac{1}{2}\int_{\R^n}\left|\nabla^{m}u\right|^2=\begin{cases}
\frac{1}{2}\int_{\R^n}|\Delta^{\frac{m}{2}}u|^2  &\mbox{if }m\mbox{ is even}\\
\frac{1}{2}\int_{\R^n}|\nabla\Delta^{\frac{m-1}{2}}u|^2 &\mbox{if }m\mbox{ is odd}
\end{cases}
\eex for any $u\in W^{m,2}(\R^n,N)$,
where $\Delta$ is the Laplace operator on $\R^n$ and
$$W^{m,2}(\R^n,N)=\left\{v\in W^{m,2}(\R^n,\R^l):
\ v(x)\in N\mbox{ for a.e. }x\in\Omega \right\}.$$ 
Recall that a map
$u\in W^{m,2}(\R^n,N)$ is called a polyharmonic map if $u$ is the
critical point of $E_m$. The Euler-Lagrange equation of polyharmonic maps is (see
Gastel-Scheven \cite{Gastel-Scheven}):
\beq\label{eq-polyharmonic}
\begin{split}(-1)^m \D^m u=F(u):=&(-1)^m {\rm{div}}^m\left(
\sum\ls_{k=0}^{m-2}\left(\begin{matrix}m-1\\k \end{matrix}\right)\nabla^{m-k-1}(\Pi(u))\nabla^{k+1}u\right)\\
&-\sum\ls_{k=0}^{m-1}(-1)^k \left(\begin{matrix}m\\k \end{matrix}\right){\rm{div}}^k\left(\nabla^{m-k}(\Pi(u))\nabla^m u\right)
\end{split}\eeq 
where $\Pi: N_{\delta}\to N$ is the nearest point projection from the $\delta$-neighborhood
of $N$ to $N$, which is smooth provide $\delta=\delta(N)>0$ is sufficiently small.
It is readily seen that (\ref{eq-polyharmonic})
becomes the equation of harmonic maps for $m=1$,  and of extrinsic biharmonic maps
for $m=2$.

Motivated by the study of heat flow of harmonic and
biharmonic maps,  we consider the heat
flow of polyharmonic maps, i.e. $u:\mathbb{R}^n\times
\mathbb{R}_+\rightarrow N$ solves

\begin{align}
u_t+(-1)^m\D^mu=&F(u) \label{heat-eq-polyharmonic}\ \ \ {\rm{in}} \ \mathbb R^n\times (0,+\infty)\\
u\big|_{t=0}=&u_0\label{initial-data} \ \ \ \ \ \  {\rm{on}}\ \mathbb R^n,
\end{align}
where $u_0:\R^n\rightarrow N$ is a given map.

The heat flow of harmonic maps, (\ref{heat-eq-polyharmonic}) for $m=1$,  
has been extensively studied.  For smooth initial data, the existence of 
global smooth heat flow of harmonic maps has been established by (i) 
Eells-Sampson \cite {Eells-Sampson} under the assumption that
the sectional curvature $K_N\le 0$, and (ii) Hildebrandt-Kaul-Widman
\cite{Hildebrandt-Kaul-Widman}  under the assumption that the image of $u_0$ 
is contained in a geodesic ball
$B_R$ in $N$ with radius $R<\frac{\pi}{2\sqrt{\max\ls_{B_R}|K_N|}}$.
In general, the short time smooth heat flow of harmonic maps may
develop singularity at finite time, see
Coron-Ghidaglia \cite{Coron-Ghidaglia}, Chen-Ding \cite{Chen-Ding},
and Chang-Ding-Ye \cite{Chang-Ding-Ye}.  However, Chen-Struwe
\cite{Chen-Struwe} (see also Chen-Lin \cite{Chen-Lin} and Lin-Wang \cite{Lin-Wang}) proved
the existence of partially smooth, global weak solutions to
(\ref{heat-eq-polyharmonic})-(\ref{initial-data}) for smooth initial
data $u_0$. For rough initial data $u_0$, the second author recently
proved in \cite{Wang-harmonic} the well-posedness for the heat flow of harmonic maps
provided the BMO norm of $u_0$ is  small.

When $m=2$, (\ref{heat-eq-polyharmonic}) becomes the heat flow of extrinsic
biharmonic maps, which was first studied by Lamm in
\cite{Lamm2001,Lamm2004,Lamm2005}. In particular, it
was proven in \cite{Lamm2001,Lamm2004,Lamm2005} that if $n = 4$ and
$\|u_0\|_{W^{2,2}(\R^4)}$ is sufficiently small, then there exists
a unique global smooth solution. For an arbitrary 
$u_0\in {W^{m,2}(\R^{2m})}$, it was later independently proved by
Wang \cite{Wang2007} (for $m=2$) and Gastel \cite{Gastel} (for $m\ge 2$)
that there exists a global weak solution to
(\ref{heat-eq-polyharmonic})-(\ref{initial-data}) that is smooth
away from finitely many singular times.  Very recently,
the second author established in \cite{Wang-biharmonic} 
the well-posedness for the heat flow of biharmonic maps for $u_0$
with small BMO norm.  

We would like to mention that there have been some works on the regularity
of polyharmonic maps for $m\ge 3$ in the critical dimensions $n=2m$. 
We refer the readers to Gastel-Scheven \cite{Gastel-Scheven}, Lamm-Wang
\cite{Lamm-Wang}, Goldstein-Strzelecki-Zatorska-Goldstein\cite{GSZ}, Moser \cite{Moser},
and Angelsberg-Pumberger \cite{Angelsberg-Pumberger}.

In this paper, we are interested in the well-posedness of the heat flow of
polyharmonic maps with rough initial data. In particular, we aim
to extend the techniques from  \cite{Wang-harmonic, Wang-biharmonic}
to establish the well-posedness of the heat flow of polyharmonic maps
(\ref{heat-eq-polyharmonic}) and (\ref{initial-data}) for 
$m\geq 3$ with $u_0$ having small BMO norm.

We remark that the techniques employed by Wang \cite{Wang-harmonic, Wang-biharmonic} were
motivated by the earlier work by Koch
and Tataru \cite{Koch-Tataru} on the global well-posedness of the
incompressible Navier-Stokes equation, and the recent work by Koch-Lamm
\cite{Koch-Lamm} on geometric flows with rough initial data.

We first recall the BMO spaces.  For $x\in\mathbb R^n$ and $r>0$, let 
$B_r(x)\subset\R^n$ be the ball with center  $x$ and radius $r$. 
For $f:\R^n\rightarrow \R$, let $f_{x,r}$ be the
average of $f$ over $B_r(x)$.

\begin{definition}\label{BMO}
For $f:\R^n\rightarrow \R$ and $R>0$, define
$${\rm{BMO}}_R(\R^n)=\left\{f:\R^n\rightarrow \R|\ \left[f\right]_{{\rm{BMO}}(\R^n)}
:=\sup\limits_{x\in\R^n,0<r\le R} r^{-n}\int_{B_r(x)}|f-f_{x,r}|<+\infty\right\}.$$
When $R=+\infty$, we simply write ${\rm{BMO}}(\R^n)$ for ${\rm{BMO}}_\infty(\R^n)$.
\end{definition}

For $0<T\leq\infty$, define the functional space $X_T$ by

\beq\label{X-space} X_T:=\left\{f:\R^n\times[0,T]\rightarrow\R^l\ |\
\left\|f\right\|_{X_T}:=\sup\ls_{0<t\leq
T}\|f\|_{L^{\infty}(\R^n)}+\left[f\right]_{X_T}\right\},\eeq where

\beq\label{X-seminorm} \left[f\right]_{X_T}=\sum\ls_{k=1}^{m}\{\sup\ls_{0<t\leq
T}t^{\frac{k}{2m}}\|\nabla^k
f\|_{\Li(\R^n)}+\sup\ls_{x\in\R^n,0<r\leq
T^{\frac{1}{2m}}} (r^{-n}\int_{P_r(x,r^{2m})}|\nabla^k
f|^{\frac{2m}{k}})^{\frac{k}{2m}}\}\eeq 
where $P_r(x,r^{2m})=B_r(x)\times[0,r^{2m}]$.
It is clear that $(X_T,\|\cdot\|_{X_T})$ is
a Banach space. When $T=+\infty$,
we simply write $X$, $\|\cdot\|_{X}$, and $[\cdot]_{X}$ 
for $X_{\infty}$,  $\|\cdot\|_{X_{\infty}}$, and  $[\cdot]_{X_{\infty}}$
respectively.

The main theorem is

\btm \label{main-theorem}{\it There exists an
$\va_0>0$ such that for any $R>0$
 if $u_0:\mathbb R^n\to N$ has $[u_0]_{{\rm{BMO}}_R(\R^n)}\leq\va_0$, then there exists a unique global
solution $u:\R^n\times [0,R^{2m}]\rightarrow N$ to
(\ref{heat-eq-polyharmonic}) and (\ref{initial-data}) with small semi-norm
$[u]_{X_{R^{2m}}}$.} \etm

As a direct consequence, we have
\begin{corollary}
There exists an $\va_0>0$ such that if $u_0:\mathbb R^n\to N$ has
$[u_0]_{{\rm{BMO}}(\R^n)}\leq\va_0$, then there exists a unique global
solution $u:\R^n\times \mathbb R_+\rightarrow N$ to
(\ref{heat-eq-polyharmonic}) and (\ref{initial-data}) with small semi-norm
$[u]_{X}$.
\end{corollary}

We follow  the arguments in \cite{Wang-harmonic, Wang-biharmonic} very closely.
The paper is written as follows. In section 2, we present some basic estimates on the polyharmonic
heat kernel. In section 3, we present some crucial estimates on the polyharmonic
heat equation. In section 4, we prove Theorem 1.2. 

\section {The polyharmonic heat kernel}
\setcounter{equation}{0}

In this section, we will prove some basic properties on the
polyharmonic heat kernel.

The fundamental solution of the polyharmonic heat equation:
\begin{equation} b_t(x,t)+(-1)^m\D^mb(x,t)=0\ \mbox{in}\ \R^n\times\R_+\end{equation}
is given by 
\begin{equation}\label{poly_kernel}
b(x,t)=t^{-\frac{n}{2m}}g\left(\frac{x}{t^{\frac{1}{2m}}}\right),
\end{equation}
where
\beq \label{poly2.1}
g(x)=(2\pi)^{-\frac{n}{2}}\int_{\R^n}e^{ix\cdot\xi-|\xi|^{2m}}d\xi,
\ x\in\mathbb R^n. \eeq
It is easy to see that $g$ is smooth, radial, and
\begin{proposition} For any $L\ge 0, k\ge 0$, there exists $C=C(k,L)>0$ such that
\begin{equation}\label{poly_decay}
|\nabla^k g(x)|\le C (1+|x|)^{-L}, \ \forall x\in\mathbb R^n.
\end{equation}
\end{proposition}
\noindent{\it Proof}. For $k\ge 0$ and $L\ge 0$, since
$$
\nabla^k_x \left(e^{ix\cdot\xi}\right)={(i\xi)^k}{(ix)^{-L}} \nabla^L_{\xi} \left(e^{ix\cdot\xi}\right),
$$
we have, by integration by parts,
\begin{eqnarray*}
|\nabla^k g(x)|
&=&\left|\int_{\mathbb R^n} (ix)^{-L} e^{ix\cdot\xi} \nabla^L_\xi\left((i\xi)^k e^{-|\xi|^{2m}}\right)\,d\xi\right|\\
&\leq& C(k,L)(1+|x|)^{-L}.
\end{eqnarray*}
This completes the proof. \qed\\

As a direct consequence of (\ref{poly_decay}), we have the following properties for the polyharmonic heat kernel $b$
\blm \label{poly-heat-kernel-lemma} For any $k,L\ge 0$, there exist $C_1>0$ depending on  $n, L$ and
$C_2,C_3>0$ depending on $n,k,L$ such that for any $x\in\R^n$ and $t>0$, it holds:

\beq\label{poly2.2} |b(x,t)|\leq
C_1 t^{-\frac{n}{2m}}\left(1+\frac{|x|}{t^{\frac1{2m}}}\right)^{-L},
\eeq

\beq\label{poly2.3} |\nabla^kb(x,t)|\leq
C_2\left(t^{-\frac{1}{2m}}\right)^{n+k-L}\left(t^{\frac{1}{2m}}+|x|\right)^{-L}, \eeq

\beq\label{poly2.4}\|\nabla^k b(x,t)\|_{L^1(\R^n)}\leq
C_3t^{-\frac{k}{2m}}.
\eeq
\elm

At the end of this section, we recall that the solution to the
Dirichlet problem of inhomogeneous polyharmonic heat equation

\begin{align}
u_t(x,t)+(-1)^m\D^mu(x,t)=&f(x,t)\ \mbox{in}\ \R^n\times\R_+,
\label{inho-eq-polyharmonic}\\
u(x,0)=&u_0(x)\ \mbox{on}\ \R^n\label{inho-initial-data}
\end{align} is given by the following Duhamel formula:

\beq\label{poly2.6} u=\mathbf{G}u_0+\mathbf{S}f, \eeq
where

\beq\label{poly2.7} {\mathbf G}u_0(x,t):=\int_{\R^n}b(x-y,t)u_0(y)dy,\
(x,t)\in\R^n\times\R_+, \eeq
and

\beq\label{poly2.8}
{\mathbf S}f(x,t):=\int_{0}^{t}\int_{\R^n}b(x-y,t-s)f(y,s)dyds,\
(x,t)\in\R^n\times\R_+. \eeq

\section {Basic estimates for the polyharmonic heat equation}
\setcounter{equation}{0}

In this section, we will provide some crucial estimates for the
solution of the polyharmonic heat equation with initial data in BMO
spaces.

\blm\label{poly-eq-lemma} {For $0<R\le +\infty$, if $u_0\in {\rm{BMO}}_R(\R^n)$, then
$\hat{u}_0:=\mathbf{G}u_0$ satisfies

\beq\label{poly3.1}\sum\ls_{k=1}^{m}\sup\ls_{x\in\R^n,0<r\le R}
r^{-n}\int_{P_r(x,r^{2m})}r^{2k-2m}|\nabla^{k}\hat{u}_0|^2\leq
C\left[u_0\right]^2_{{\rm{BMO}}_R(\R^n)},\eeq
and
\beq\label{poly3.2} \sum\ls_{k=1}^{m}\sup\ls_{0<t\le R^{2m}}t^{\frac{k}{2m}}
\left\|\nabla^k\hat{u}_0(t)\right\|_{L^{\infty}(\R^n)}\leq
C\left[u_0\right]_{{\rm{BMO}}_R(\R^n)}.\eeq
If, in addition, $u_0\in L^\infty(\R^n)$ then
\beq\label{poly3.3}\sum\ls_{k=1}^{m-1}\sup\ls_{x\in\R^n,0<r\le R}
r^{-n}\int_{P_r(x,r^{2m})}|\nabla^{k}\hat{u}_0|^{\frac{2m}{k}}\leq
C\left\|u_0\right\|^{\frac{2m}{k}-2}_{L^{\infty}(\R^n)}\cdot\left[u_0\right]^2_{{\rm{BMO}}_R(\R^n)},\eeq
}\elm

The proof of Lemma 3.1 is similar to \cite{Wang-biharmonic} Lemma 3.1. 
For completeness, we sketch it here. Let $\mathcal{S}$ denote the class of Schwartz functions,
the following characterization of BMO spaces, due to Carleson,  is well-known (see, Stein \cite{Stein}).

\blm\label{Carleson-BMO}{\it For $0<R\le+\infty$, 
let $\Phi\in \mathcal{S}$ be such that $\int_{\R^n}\Phi=0$ and denote for $t>0$,
$\P_t(x)=t^{-n}\P(\frac{x}{t})$, $x\in\R^n$. If $f\in {\rm{BMO}}_R(\R^n)$,
then 
\beq\label{Carleson-BMO-inq}
\sup\ls_{x\in\R^n,0<r\le R}r^{-n}\int_0^r\int_{B_r(x)}
|\P_t\ast f|^2(x,t)\frac{dxdt}{t}\leq C\left[u_0\right]_{{\rm{BMO}}_R(\R^n)}
\eeq
for some $C=C(n)>0.$}\elm

\noindent{\it Proof of Lemma 3.1}.  Let $g$ be given by (\ref{poly2.1}) and $\P^i=\nabla^ig$ for
$i=1,\cdots,m$. Then $\P^i\in\mathcal{S}$ and
$\int_{\R^n}\P^i=0$ for  $i=1,\cdots,m$. Direct calculations show
$$\P^i_t(x)=t^{-n}(\nabla^ig)\left(\frac{x}{t}\right)
=t^i\nabla^i\left(t^{-n}g(\frac{x}{t})\right)=t^i\nabla^ig_t(x),$$
where $g_t(x)=t^{-n}g(\frac{x}{t})$.
Hence we have
$$\P^i_t\ast u_0(x)=t^i\nabla^i\left(g_t\ast u_0\right)(x).$$
Since the polyharmonic heat kernel $b(x,t)=g_{t^{\frac{1}{2m}}}(x)$, we have
$$\P^i_t\ast u_0(x)=t^i\nabla^i[(b(\cdot,t^{2m})\ast u_0)(x)]
=t^i\nabla^i(\mathbf Gu_0)(x,t^{2m}).$$
Thus Lemma 3.1 implies that for $i=1,\cdots,m$,
\bex
\begin{split}
C\left[u_0\right]_{{\rm{BMO}}_R(\R^n)}^2\geq&\sup\ls_{x\in\R^n,0<r\le R}
r^{-n}\int_{0}^{r}\int_{B_r(x)}|\P^i_t\ast u_0|^2\frac{dxdt}{t}\\
=&\sup\ls_{x\in\R^n,0<r\le R}
r^{-n}\int_{0}^{r}\int_{B_r(x)}t^{2i-1}|\nabla^i \mathbf G u_0|^2(x,t^{2m})dxdt\\
=&\frac{1}{2m}\sup\ls_{x\in\R^n,0<r\le R}
r^{-n}\int_{P_r(x,r^{2m})}t^{\frac{2i-2m}{2m}}|\nabla^i \mathbf G u_0|^2(x,t)dxdt.
\end{split}\eex
This clearly implies (\ref{poly3.1}), since for $i=1,\cdots,m$,
$t^{\frac{2i-2m}{2m}}\geq r^{2i-2m}$ when $0<t\leq r^{2m}$.

Since $\hat{u}_0$ solves the polyharmonic heat equation:
$$({\partial_t}+(-1)^m\D^m)\hat{u}_0=0
\ {\rm{on}}\ \R^n\times(0,+\infty),$$ 
the standard theory implies
that for any $x\in\R^n$ and $r>0$,
$$\sum\ls_{k=1}^{m}r^{\frac{k}{m}}|\nabla^k\hat{u}_0|^2(x,r^{2m})
\leq C\sum\ls_{k=1}^{m}
r^{-n}\int_{P_r(x,r^{2m})}r^{2k-2m}|\nabla^{k}\hat{u}_0|^2.$$
Taking supremum over $x\in\R^n$ and $0<t=r^{2m}\le R^{2m}$ yield (\ref{poly3.2}).

For (\ref{poly3.3}), observe that $u_0\in\Li(\R^n)$ implies 
$\P^{i}_t\ast u_0\in\Li(\R^n)$ for $i=1,\cdots,m-1$, and
$$\|\P^{i}_t\ast u_0\|_{\Li(\R^n)}\leq\|\P^{i}\|_{L^{1}(\R^n)}\|u_0\|_{\Li(\R^n)}
\leq \|\nabla^ig\|_{L^{1}(\R^n)}\|u_0\|_{\Li(\R^n)}
\leq C\|u_0\|_{\Li(\R^n)}.$$
Hence
\bex\begin{split}
&\sup\ls_{x\in\R^n,0<r\le R}
r^{-n}\int_{P_r(x, r^{2m})}|\nabla^{i}\mathbf Gu_0|^{\frac{2m}{i}}dxdt\\
&=\sup\ls_{x\in\R^n,0<r\le R}
r^{-n}\int_{B_r(x)\times[0,r]}|\P_t^i\ast u_0|^{\frac{2m}{i}}\frac{dxdt}{t}\\
&\leq\left(\sup\ls_{t>0}\|\P^{i}_t\ast u_0\|_{\Li(\R^n)}
\right)^{\frac{2m}{i}-2}_{L^{\infty}(\R^n)}
\cdot\sup\ls_{x\in\R^n,0<r\le R}
r^{-n}\int_{B_r(x)\times [0,r]}|\P_t^i\ast u_0|^{2}\frac{dxdt}{t}\\
&\leq
C\left\|u_0\right\|^{\frac{2m}{i}-2}_{L^{\infty}(\R^n)}\cdot\left[u_0\right]^2_{{\rm{BMO}}_R(\R^n)}
\end{split}\eex
This implies (\ref{poly3.3}).
\endpf\\

Now we prove an important estimate on the distance of $\hat{u}_0$ to the
manifold $N$ in term of the BMO norms of $u_0$. More
precisely,

\blm\label{poly-tangent-lemma} For any $\delta>0$, there exists
$K_0=K_0(\delta,N)>0$ such that for $0<R\le+\infty$, if
$u_0\in{\rm{BMO}}_R(\R^n)$ then 

\beq\label{poly-tangent-inq}\mbox{dist}(\hat{u}_0(x,t),N)\leq
K_0[u_0]_{{\rm{BMO}}_R(\R^n)}+\delta,\ \forall x\in\R^n,  \ 0\le t\le (\frac{R}{K_0})^{2m}.
\eeq
\elm

\noindent{\it Proof}.
 For any $x\in\R^n$, $t>0$ and $K>0$, denote
$$c^{K}_{x,t}=\frac{1}{|B_K(0)|}\int_{B_K(0)}u_0(x-t^{\frac{1}{2m}}z)dz.$$
Let $g$ be given by (\ref{poly2.1}). Then, by a change of variables, we have
$$\hat{u}_0(x,t)=\int_{\R^n}g(y)u_0(x-t^{\frac{1}{2m}}y)dy.$$
Applying (2.5) (with $L=n+1$) from Lemma 2.2, we have
\beq\label{poly3.4}\begin{split}
\left|\hat{u}_0(x,t)-c_{x,t}^{K}\right|\leq&
\int_{\R^n}g(y)\left|u_0(x-t^{\frac{1}{2m}}y)-c_{x,t}^{K}\right|dy\\
\leq&\left\{\int_{B_K(0)}+\int_{\R^n\setminus B_K(0)}\right\}
g(y)\left|u_0(x-t^{\frac{1}{2m}}y)-c_{x,t}^{K}\right|dy\\
\leq&\int_{B_K(0)}
\left|u_0(x-t^{\frac{1}{2m}}y)-c_{x,t}^{K}\right|dy\\
&+C\|u_0\|_{\Li(\R^n)}\int_{\R^n\setminus B_K(0)}\frac{1}{|y|^{n+1}}\,dy\\
\leq&K^n\left[u_0\right]_{{\rm{BMO}}_{Kt^{\frac{1}{2m}}}(\R^n)}+\delta,
\end{split}\eeq
provided we choose a sufficiently large $K=K_0(\delta,N)>0$ so that
$$C\|u_0\|_{\Li(\R^n)}\int_{\R^n\setminus B_K(0)}
\frac{1}{|y|^{n+1}}\,dy\leq \delta.$$
On the other hand, since $u_0(\R^n)\subset N$, we have
\beq\label{poly3.5}{\rm{dist}}(c^K_{x,t},N)\leq \frac{1}{|B_K(0)|}
\int_{B_K(0)}\left|c^K_{x,t}-u_0(x-t^{\frac{1}{2m}}y)\right|dy
\leq \left[u_0\right]_{{\rm{BMO}}_{Kt^{\frac1{2m}}}(\R^n)}.
\eeq
Putting (\ref{poly3.4}) and (\ref{poly3.5}) together yields (\ref{poly-tangent-inq})
holds for $t\le (\frac{R}{K})^{2m}$.
\endpf

\section {\bf Boundedness of the operator $\mathbf{S}$}
\setcounter{equation}{0}

In this section, we introduce several function spaces and establish
the boundedness of the operator $\mathbf S$ between these spaces.

For $0<T<\infty$, the spaces $Y_T^k$, for $k=0,\cdots,m-1$, are the sets consisting of
all functions $f:\R^n\times[0,T]\rightarrow \R$ such that

\beq\label{Yk-norm} \|f\|_{Y_T^{k}}:=\sup\ls_{0<t\leq
T}t^{\frac{2m-k}{2m}}\|f\|_{\Li(\R^n)} +\sup\ls_{x\in\R^n,0<R\leq
T^{\frac{1}{2m}}}
\left(R^{-n}\int_{P_R(x)}|f|^{\frac{2m}{2m-k}}\right)^{\frac{2m-k}{2m}}.
\eeq 
Notice that  $(Y^k_T,\|\cdot\|_{Y_T^k})$ is a Banach
space for $k=0,\cdots,m-1$. When $T=+\infty$, we simply denote $(Y^k,\|\cdot\|_{Y^k})$
for $(Y_{\infty}^{k},\|\cdot\|_{Y^k_\infty})$.

Let the operator $\mathbf S$ be defined by (\ref{poly2.8}). Then we have

\blm\label{bounded-of-S}{\it For any $0<T\le +\infty$ and $k=0,\cdots,m-1$, 
if $f\in Y^k_T$, then $\mathbf S(\nabla^{\alpha}f)\in X_T$ and

\beq\label{bounded-inq-S}
\left\|\mathbf S(\nabla^{\alpha}f)\right\|_{X_T}\leq C\left\|f\right\|_{Y^k_T},
\eeq
where $\alpha=(\alpha_1,\cdots\alpha_n)$ is any multi-index of order  $k$.

}\elm

\vspace{2mm}

\pf We need to show the point wise estimate

\beq\label{poly4.1} \sum\ls_{i=0}^{m}R^i|\nabla^i
\mathbf S(\nabla^{\alpha}f)|(x,R^{2m})\leq C\|f\|_{Y^k_T},\quad\forall
x\in\R^n,\ 0<R\le T^{\frac{1}{2m}}, \eeq 
and the integral estimate for $0<R\le T^{\frac{1}{2m}}$:

\beq\label{poly4.2} \sum\ls_{i=1}^{m}R^{-\frac{in}{2m}}\left\|\nabla^i
\mathbf S(\nabla^{\alpha}f)\right\|_{L^{\frac{2m}{i}}(P_R(x,R^{2m}))}\leq C\left\|f\right\|_{Y^k_T}
\eeq
By suitable scaling, we may assume $T\ge 1$.
Since both estimates are translation and scale invariant, it
suffices to show (\ref{poly4.1}) and (\ref{poly4.2}) hold for $x=0$
and $R=1$.

For $i=0,\cdots,m$ and $\alpha=(\alpha_1,\cdots,\alpha_n)$ with order $k$,
we have

\bex\begin{split} \left|\nabla^i\mathbf S(\nabla^{\alpha}f)\right|(0,1)
=&\left|\int_0^1\int_{\R^n}\nabla^{i+\alpha}b(y,1-s)f(y,s)dyds\right|\\
\leq&\left\{\int_{\frac12}^1\int_{\R^n}+\int^{\frac{1}{2}}_{0}\int_{B_2}
+\int^{\frac{1}{2}}_{0}\int_{\R^n\setminus
B_2}\right\}\left|\nabla^{i+k}b(y,1-s)\right|\left|f(y,s)\right|dyds\\
=&I_1+I_2+I_3.
\end{split}\eex
Applying Lemma \ref{poly-heat-kernel-lemma}, we can estimate $I_1$,
$I_2$ and $I_3$ as follows.

\bex\begin{split} |I_1|\leq&\left(\sup\ls_{\frac{1}{2}\leq s\leq
1}\|f(s)\|_{\Li(\R^n)}\right)\left(\int_{\frac{1}{2}}^{1}\left\|\nabla
^{i+k}b(\cdot,1-s)\right\|_{L^1(\R^n)}ds\right)\\
\leq& C\|f\|_{Y^k_1}\int^{\frac{1}{2}}_{0}s^{-\frac{i+k}{2m}}ds
\ \ ({\rm{by}}\ (2.7))\\
\leq&C\|f\|_{Y^k_1}\ \ \ ({\rm{since}}\ i+k\le 2m-1).
\end{split}\eex

\bex\begin{split} |I_2|\leq&\left(\sup\ls_{0\leq s\leq
\frac{1}{2}}\|\nabla^{i+k}b(\cdot,1-s)\|_{\Li(\R^n)}\right)
\left(\int_{B_2\times[0,\frac{1}{2}]}|f(y,s)|dyds\right)\\
\leq& C\int_{B_2\times[0,\frac{1}{2}]}|f(y,s)|dyds\\
\leq&C\|f\|_{Y^k_1}.
\end{split}\eex

\bex\begin{split} |I_3|\leq&\int_{0}^{\frac{1}{2}}\int_{\R^n\setminus
B_2}\left|\nabla^{i+k}b(y,1-s)\right||f(y,s)|dyds\\
\leq& C\int_{0}^{\frac{1}{2}}\int_{\R^n\setminus
B_2}|y|^{-(n+1)}|f(y,s)|dyds\ \ \  ({\rm{by}} \ (2.6) \ {\rm{for}}\ L=n+1)\\
\leq&\left(\sum\ls_{k=2}^{\infty}k^{n-1} k^{-(n+1)}\right)
\left(\sup\ls_{x\in\R^n}\int_{P_1(x,1)}|f(y,s)|dyds\right)\\
\leq&C\left(\sum_{k=2}^\infty k^{-2}\right)\|f\|_{Y^k_1}
\leq \|f\|_{Y^k_1}.
\end{split}\eex

Now we want to show (\ref{poly4.2}) by the energy method. Denote
$w=\mathbf S(\nabla^{\alpha}f)$. Then $w$ solves

\beq\label{poly4.3} ({\partial_t}+(-1)^m\D^m)w=\nabla^{\alpha}f\quad\mbox{in
}\R^n\times(0,+\infty);\quad w|_{t=0}=0.
\eeq
Let $\eta\in
C^{\infty}_{0}(B_2)$ be a cut-off function of $B_1$. Multiplying
(\ref{poly4.3}) by $\eta^4w$ and integrating over $\R^n\times[0,1]$,
we obtain

\beq\label{poly4.4} \int_{\R\times\{1\}}|w|^2\eta^4
+2\int_{\R^n\times[0,1]}\nabla^mw\cdot\nabla^m(w\eta^4)
=2\int_{\R^n\times[0,1]}\nabla^{\alpha}f\cdot w\eta^4. \eeq
By the H\"older inequality, we have
\beq\label{poly4.5}\begin{split}
&\int_{\R^n\times[0,1]}\nabla^mw\cdot\nabla^m(w\eta^4)\\
=&\int_{\R^n\times[0,1]}|\nabla^m(w\eta^2)|^2
+\int_{\R^n\times[0,1]}\nabla^mw\left(\sum\ls_{\beta=0}^{m-1}\nabla^{\beta}(w\eta^2)
\cdot\nabla^{m-\beta}(\eta^2)\right)\\
&-\int_{\R^n\times[0,1]}\nabla^m(w\eta^2)\left(\sum\ls_{\beta=0}^{m-1}\nabla^{\beta}(w)
\cdot
\nabla^{m-\beta}(\eta^2)\right)\\
\geq&\frac{1}{2}\int_{\R^n\times[0,1]}|\nabla^m(w\eta^2)|^2
-C\sum\ls_{\beta=0}^{m-1}\int_{B_2\times[0,1]}|\nabla^{\beta}w|^2
\end{split}\eeq

\beq\label{poly4.6}\begin{split}
&\int_{\R^n\times[0,1]}\nabla^{\alpha}f\cdot w\eta^4
=(-1)^k\int_{\R^n\times[0,1]}f\cdot \nabla^{\alpha}[(w\eta^2)\eta^2]\\
\leq&C\sum\ls_{\beta=0}^{k}\int_{\mathbb R^n\times[0,1]}|f|
|\nabla^{\beta}( w\eta^2)|\\
\leq&C\sum\ls_{\beta=0}^{k-1}
\sup\ls_{0<t\leq 1}t^{\frac{2m-k}{2m}}\|f\|_{\Li(\R^n)}
\cdot\sup\ls_{0<t\leq 1}t^{\frac{\beta}{2m}}\|\nabla^{\beta}w\|_{\Li(\R^n)}
\cdot\int_0^1 t^{-1+\frac{k-\beta}{2m}}dt\\
&+C\|f\|_{L^{\frac{2m}{2m-k}}(B_2\times[0,1])}\cdot
\|\nabla^{k}( w\eta^2)\|_{L^{\frac{2m}{k}}(\R^n\times[0,1])}\\
\leq&C\|f\|_{Y^k_1}^2+C\|f\|_{Y^k_1}\cdot \|\nabla^{k}(
w\eta^2)\|_{L^{\frac{2m}{k}}(\R^n\times[0,1])}
\end{split}\eeq
To estimate the last term, we need the Nirenberg interpolation
inequality: for $k\le m-1$,

$$\|\nabla^{k}(w\eta^2)\|_{L^{\frac{2m}{k}}(\R^n)}^{\frac{2m}{k}}\leq
C\|w\eta^2\|_{\Li(\R^n)}^{\frac{2m}{k}-2}
\|\nabla^{m}(w\eta^2)\|_{L^{2}(\R^n)}^{2},$$ which, after
integrating with respect to $t\in[0,1]$, implies

\beq\label{poly5.13}
\|\nabla^{k}(w\eta^2)\|_{L^{\frac{2m}{k}}(\R^n\times[0,1])}\leq
C\sup\ls_{0\leq t\leq 1}\|w\|_{\Li(\R^n)}^{1-\frac{k}{m}}
\|\nabla^{m}(w\eta^2)\|_{L^{2}(\R^n\times[0,1])}^{\frac{k}{m}},\eeq
Putting (\ref{poly5.13}), (\ref{poly4.5}) and (\ref{poly4.6}) into
(\ref{poly4.4}), we have

\beq\label{poly4.7}
\begin{split}
&\int_{\R^n\times[0,1]}|\nabla^m(w\eta^2)|^2\\
\leq&
C\sum\ls_{\beta=0}^{m-1}\int_{B_2\times[0,1]}|\nabla^{\beta}w|^2
+C\|f\|_{Y^k_1}^2+C\|f\|_{Y^k_1}\cdot \|\nabla^{k}(
w\eta^2)\|_{L^{\frac{2m}{k}}(\R^n\times[0,1])}\\
\leq&C\sum\ls_{\beta=0}^{m-1}\left[\int_0^1t^{-\frac{\beta}{m}}dt
\cdot\sup\ls_{0<t\leq1}(t^{\frac{\beta}{m}}\|\nabla^{\beta}w(t)\|_{\Li(\R^n)}^2)\right]\\
&+C\|f\|_{Y^k_1}^2+C\|f\|_{Y^k_1}\sup\ls_{0\leq t\leq
1}\|w\|_{\Li(\R^n)}^{1-\frac{k}{m}}
\|\nabla^{m}(w\eta^2)\|_{L^{2}(\R^n\times[0,1])}^{\frac{k}{m}}\\
\leq&C\|f\|_{Y^k_1}^2+\frac{1}{2}\int_{\R^n\times[0,1]}|\nabla^m(w\eta^2)|^2
+C\|f\|_{Y^k_1}^q\cdot \|w\|_{\Li(\R^n)}^{(1-\frac{k}{m})q}\\
\leq&\frac{1}{2}\int_{\R^n\times[0,1]}|\nabla^m(w\eta^2)|^2+C\|f\|_{Y^k_1}^2,
\end{split}\eeq where $q=\frac{2m}{2m-k}$.
Therefore, we obtain

\beq\label{poly4.8} \int_{P_1(0,1)}|\nabla^mw|^2
\leq\int_{\R^n\times[0,1]}|\nabla^m(w\eta^2)|^2\leq C\|f\|_{Y^k_1}^2.
\eeq 
For  $i=1,\cdots,m-1$, applying Nirenberg's interpolation inequality gives

\beq\label{poly4.9}\begin{split} &\int_{P_1(0,1)}|\nabla^i
w|^{\frac{2m}{i}}\leq\int_{\R^n\times[0,1]}|\nabla^i
(w\eta^2)|^{\frac{2m}{i}}\\
\leq& \sup\ls_{0\leq t\leq 1}\|w\|_{\Li(\R^n)}^{\frac{2m}{i}-2}
\|\nabla^{m}(w\eta^2)\|_{L^{2}(\R^n\times[0,1])}^{2}\\
\leq&C\|f\|_{Y^k_1}^{\frac{2m}{i}-2}\cdot\int_{\R^n\times[0,1]}|\nabla^m
(w\eta^2)|^{2}\\
\leq&C\|f\|_{Y^k_1}^{\frac{2m}{i}}
\end{split}\eeq
(\ref{poly4.8}) and (\ref{poly4.9}) imply (\ref{poly4.2}). This
completes the proof.
\endpf

\section {\bf Proof of Theorem 1.2}
\setcounter{equation}{0}

This section is devoted to the proof of Theorem \ref{main-theorem}.
The idea is based on the fixed point theorem in a small ball inside
$X_{R^{2m}}$.

Since the image of a map $u\in X_{R^{2m}}$ may not be contained in $N$,
we first need to extend $\Pi$ to $\mathbb R^l$, denoted as $\widetilde{\Pi}$, such that
$\widetilde{\Pi}\in C^\infty(\mathbb R^l)$ and $\widetilde{\Pi}\equiv\Pi$ in $N_{\delta_N}$.

Let
\beq\label{poly5.9} \begin{split} \widetilde{F}(u):
=&(-1)^m{\rm{div}}^m\left(\sum\ls_{k=0}^{m-2}\left(\begin{matrix}m-1\\k \end{matrix}\right) \nabla^{m-k-1}(\widetilde{\Pi}(u))\nabla^{k+1} u\right)\\
&-\sum\ls_{k=0}^{m-1}(-1)^k\left(\begin{matrix}m\\k \end{matrix}\right){\rm{div}}^k\left(\nabla^{m-k}(\widetilde{\Pi}(u))\nabla^mu\right).
\end{split}\eeq

For $k=0,\cdots,m-2$, define
$$F_{k}(u)=(-1)^{k+1}\left(\begin{matrix}m\\k \end{matrix}\right)\nabla^{m-k}(\widetilde{\Pi}(u))\nabla^mu.
$$
For $k=m-1$, define
\begin{eqnarray*}&&F_{m-1}(u)\\
&=&(-1)^m \left(\begin{matrix}m\\m-1 \end{matrix}\right)\nabla(\widetilde{\Pi}(u))\nabla^mu\\
&+&(-1)^m \sum_{k=0}^{m-2}\left(\begin{matrix}m-1\\k \end{matrix}\right)
\left({\rm{div}}(\nabla^{m-k-1}(\widetilde{\Pi}(u)))\nabla^{k+1}u+
\nabla^{m-k-1}(\widetilde{\Pi}(u)){\rm{div}}(\nabla^{k+1}u)\right).
\end{eqnarray*}
Then (\ref{heat-eq-polyharmonic}) can be written
as
\beq\label{rewrite-eq-poly} ({\partial_t}+(-1)^m\D^m)u=\sum\ls_{k=0}^{m-1}{\rm{div}}^k\left(F_{k}(u)\right). \eeq 
The first observation is

\blm\label{poly-lemma5.1}{\it For $0<R\le+\infty$, if $u\in X_{R^{2m}}$, 
then $F_k(u)\in Y^k_{R^{2m}}$ for $k=0,\cdots,m-1$. Moreover, there exists
$C>0$ depends on $N$ and $\|u\|_{L^\infty(\mathbb R^n\times [0,R^{2m}])}$ such that

\beq\label{poly5.1} \left\|F_k(u)\right\|_{Y^k_{R^{2m}}}\leq C
\sum_{l=1}^m \left[u\right]_{X_{R^{2m}}}^{\frac{2m-k}{l}},
\ 0\le k\le m-1. \eeq

}\elm

\noindent{\it Proof}. It follows from  direct calculations and H\"older
inequality that
$$|F_k(u)|\le C\sum_{l=1}^m |\nabla^l u|^{\frac{2m-k}{l}},
\ 0\le k\le m-1,$$
where $C>0$ depends on $N$ and $\|u\|_{L^\infty(\mathbb R^n\times [0,R^{2m}])}$.
Thus we have, by the definitions of $X_{R^{2m}}$ and $Y^k_{R^{2m}}$, that 
$$\|F_k(u)\|_{Y^k_{R^{2m}}}
\le C\sum_{l=1}^m \left[u\right]_{X_{R^{2m}}}^{\frac{2m-k}{l}},
\ \forall 0\le k\le m-1.$$
This completes the proof. \endpf\\

By the Duhamel formula (\ref{poly2.6}), the solution $u$ to
(\ref{rewrite-eq-poly}) is given by

\beq\label{poly5.2}
u(x,t)=\mathbf Gu_0+\sum\ls_{k=0}^{m-1}\mathbf S\left({\rm{div}}^k(F_k(u))\right). \eeq
From now on,denote 
$$\hat{u}_0=\mathbf Gu_0.$$ 
Define the mapping operator $\T$ on $X_{R^{2m}}$
by letting

\beq\label{poly5.3} \T
u(x,t)=\hat{u}_0+\sum\ls_{k=0}^{m-1}\mathbf S\left({\rm{div}}^k(F_k(u))\right). \eeq

The following property follows directly from Lemma
\ref{poly-eq-lemma}.

\blm\label{poly-lemma5.2}{\it For any $0<R\le +\infty$,
if $u_0:\R^n\rightarrow N$, then $\hat{u}_0\in X_{R^{2m}}$, and 

\beq\label{poly5.4} \|\hat{u}_0\|_{\Li(\R^{n}\times [0, R^{2m}])}\leq C
\|u_0\|_{\Li(\R^{n})},\quad \left[\hat{u}_0\right]_{X_{R^{2m}}}\leq
C\left[u_0\right]_{{\rm{BMO}}_R(\R^n)}.\eeq
}\elm \vspace{2mm}

For any $\va>0$ and $0<R\le +\infty$, let

$$\B_{\va}(\hat{u}_0):=\left\{u\in X:\  \left\|u-\u\right\|_{X_{R^{2m}}}\leq \va\right\}$$ 
be the ball in $X_{R^{2m}}$ with center $\u$ and radius $\va$. 
By the triangle inequality, we have

\beq\label{poly5.5} \|u\|_{\Li(\R^{n}\times [0,R^{2m}])}\leq C
\|u_0\|_{\Li(\R^{n})}+\va,\ \left[u\right]_{X_{R^{2m}}}\leq
C\left[u_0\right]_{{\rm{BMO}}_R(\R^n)}+\va,\ \forall u\in\B_{\va}(\u). \eeq 
Thus we have

\blm\label{poly-lemma5.3}{\it For $0<R\le +\infty$, if $u_0:\R^n\rightarrow N$ has
$[u_0]_{{\rm{BMO}}_R(\R^n)}\leq\va$, then

\beq\label{poly5.6} \|u\|_{\Li(\R^{n+1}_{+})}\leq C +\va,\quad
[u]_{X}\leq C\va,\quad\forall u\in\B_{\va}(\u),\eeq for some
$C=C(n,N)>0.$

}\elm \vspace{2mm}

 The proof of Theorem \ref{main-theorem} is based
on the following two lemmas.

\blm\label{poly-lemma5.4}{\it There exists $\va_1>0$  such that 
for $0<R\le+\infty$, if $u_0:\R^n\rightarrow N$ has

$$[u_0]_{{\rm{BMO}}_R(\R^n)}\leq\va_1,$$ then 
$\T$ maps $\B_{\va_1}(\u)$ to $\B_{\va_1}(\u)$.

}\elm\vspace{2mm}

\noindent{\it Proof}. By (\ref{poly5.3}), we have
$$\mathbf T(u)-\u=\sum\ls_{k=0}^{m-1}\mathbf S\left({\rm{div}}^k (F_k(u))\right),\quad u\in\B_{\va_1}(\u).$$
Hence Lemma \ref{bounded-of-S}, Lemma \ref{poly-lemma5.1} and Lemma
\ref{poly-lemma5.2} imply that for any $u\in\B_{\va_1}(\u)$,

\bex \begin{split} \left\|\T(u)-\u\right\|_{X_{R^{2m}}}
&\lesssim \sum\ls_{k=0}^{m-1}\left\|\mathbf S({\rm{div}}^k (F_k(u)))\right\|_{X_{R^{2m}}}\\
&\lesssim \sum\ls_{k=0}^{m-1}\left\|F_k(u)\right\|_{Y^k_{R^{2m}}}\\
&\lesssim\sum_{k=0}^{m-1}\sum_{l=1}^m \left[u\right]_{X_{R^{2m}}}^{\frac{2m-k}{l}}\\
&\leq C\left[u\right]_{X_{R^{2m}}}^{\frac{m+1}m}\leq\va_1,
\end{split}\eex
provide $\va_1>0$ is chosen to be sufficiently small.  Hence
$\mathbf Tu\in\B_{\va_1}(\u)$. This completes the proof.
\endpf

\blm\label{poly-lemma5.5}{\it There exist $0<\va_2\leq\va_1$ and
$\theta_0\in(0,1)$ such that for $0<R\le+\infty$, if
$u_0:\R^n\rightarrow N$ satisfies

$$[u_0]_{{\rm{BMO}}_R(\R^n)}\leq\va_2,$$ then
$\T:\B_{\va_2}(\u):\rightarrow\B_{\va_2}(\u)$ is a
$\theta_0$-contraction map, i.e.

$$\left\|\T(u)-\T(v)\right\|_{X_{R^{2m}}}\leq\theta_0
\left\|u-v\right\|_{X_{R^{2m}}},\quad\forall u,v\in\B_{\va_2}(\u).$$

}\elm\vspace{2mm}

\noindent{\it Proof}.
 For $u,v\in\B_{\va_2}(\u)$, we have

\beq\label{poly5.10}\begin{split}
\left\|\T(u)-\T(v)\right\|_{X_{R^{2m}}}\leq&
\sum\ls_{k=0}^{m-1}\left\|\mathbf S({\rm{div}}^k(F_k(u)-F_k(v)))\right\|_{X_{R^{2m}}}\\
\lesssim&\sum\ls_{k=0}^{m-1}\left\|F_k(u)-F_k(v)\right\|_{Y^{k}_{R^{2m}}}.
\end{split}\eeq
Notice that $\|u\|_{L^\infty(\mathbb R^n\times [0,R^{2m}])}
+\|v\|_{L^\infty(\mathbb R^n\times [0,R^{2m}])}\le C_0.$
For any $k=0,\cdots,m-2$, it follows from the definition of $F_k(u)$ we have

\bex\begin{split}&|F_k(u)-F_k(v)|\\
\lesssim&
|\nabla^m(u-v)|\left[|\nabla^{m-k}u|+\sum_{j=1}^{m-k}
\left(\sum_{|\alpha|=j}(\Pi_{i=1}^n |\nabla^{\alpha_i} u|)\right)|\nabla^{m-k-j}u|\right]\\
+&|\nabla^m v|
|u-v|\left[\sum_{|\alpha|=m-k}\left(\Pi_{i=1}^n|\nabla^{\alpha_i}u|
+\Pi_{i=1}^n |\nabla^{\alpha_i}v|\right)\right]\\
+&|\nabla^m v|\left[\sum_{j=1}^{m-k}|\nabla^j(u-v)|
\left(\sum_{|\alpha|=m-k-j}(\Pi_{i=1}^n|\nabla^{\alpha_i}u|+\Pi_{i=1}^n|\nabla^{\alpha_i}v|)\right)\right].
\end{split}\eex
Hence we have

\beq\label{poly5.11}\begin{split}
\left\|F_k(u)-F_k(v)\right\|_{Y^{k}_{R^{2m}}}\lesssim&
\left[\sum_{j=1}^{m-k}([u]_{X_{R^{2m}}}^{j}+[v]_{X_{R^{2m}}}^j)\right]\left\|u-v\right\|_{X_{R^{2m}}}\\
\leq& C\va_2\left\|u-v\right\|_{X_{R^{2m}}},
\end{split}\eeq
where we have used Lemma \ref{poly-lemma5.3} in the last step.

For $k=m-1$, since
\begin{eqnarray*}
|F_{m-1}(u)-F_{m-1}(v)
&\lesssim& \left [|\nabla u||\nabla^m(u-v)|+|\nabla(u-v)||\nabla^m v|\right]\\
&+& \sum_{k=0}^{m-2}\left[|\nabla^{m-k}(\widetilde{\Pi}(u))|+|\nabla^{m-k}(\widetilde{\Pi}(v))|\right]\left|\nabla^{k+1}(u-v)\right|\\
&+&\sum_{k=0}^{m-2} \left|\nabla^{m-k}(\widetilde{\Pi}(u)-\widetilde{\Pi}(v))\right|
\left[|\nabla^{k+1}u|+|\nabla^{k+1}v|\right],
\end{eqnarray*}
we also have
\beq\label{poly5.12}\left\|F_{m-1}(u)-F_{m-1}(v)\right\|_{Y^{m-1}_{R^{2m}}}
\leq C\va_2\left\|u-v\right\|_{X_{R^{2m}}}.
\eeq Putting (\ref{poly5.11}) and (\ref{poly5.12}) into
(\ref{poly5.10}), we obtain

$$\left\|\T(u)-\T(v)\right\|_{X_{R^{2m}}}\leq C\va_2\left\|u-v\right\|_{X_{R^{2m}}}
\leq\theta_0\left\|u-v\right\|_{X_{R^{2m}}}$$
for some $\theta_0=\theta_0(\va_2)\in(0,1)$, provided $\va_2>0$ is
chosen to be sufficiently small. This completes the proof of Lemma
\ref{poly-lemma5.5}.
\endpf\\

\vspace{2mm}

\noindent{\bf Proof of Theorem \ref{main-theorem}}. It follows from
Lemma \ref{poly-lemma5.4} and Lemma \ref{poly-lemma5.5}, and the
fixed point theorem that there exists $\va_0>$ such that for
$0<R\le+\infty$,  if $[u_0]_{{\rm{BMO}}_R(\R^n)}\leq\va_0,$ then 
there exists a unique $u\in X_{R^{2m}}$
such that
$$u=\u+\mathbf S(\widetilde{F}(u))\quad \mbox{on }\R^n\times [0,R^{2m}],$$
or equivalently
$$u_t+(-1)^m\D^mu=\widetilde{F}(u)\quad \mbox{on }\R^n\times(0,R^{2m}];\ u\big|_{t=0}=u_0.$$

We want to show $u(\R^n\times [0,R^{2m}])\subset N$. By  Lemma \ref{poly-tangent-lemma},
we have that for any $x\in\mathbb R^n$ and $0\le t\le (\frac{R}{K_0})^{2m}$,

\bex\begin{split} {\rm{dist}}(u(x,t),N)\leq&{\rm{dist}}(\u,N)+\left\|u-\u\right\|_{\Li(\R^n\times [0,R^{2m}])}\\
\leq&\delta+K_0[u_0]_{{\rm{BMO}}_R(\R^n)}+\va_0\\
\leq&\delta+(1+K_0)\va_0\leq\delta_N,
\end{split}\eex
provided $\delta\leq\frac{\delta_N}{2}$ and
$\va_0=\frac{\delta_N}{2(1+K_0)}$. This yields
$u(\R^n\times [0,(\frac{R}{K_0})^{2m}])\subset N_{\delta_N}$.
Hence 
$$\widetilde{\Pi}(u)=\Pi(u), \quad \widetilde{F}(u)=F(u)
\ {\rm{on}}\ \mathbb R^n\times [0, (\frac{R}{K_0})^{2m}].
$$
Set $Q(u)=y-\Pi(y)$ for $y\in N_{\delta_N}$, and
$\rho(u)=\frac{1}{2}|Q(u)|^2$. Then direct calculations imply that
for any $y\in N_{\delta_N}$,
\begin{eqnarray*} \nabla Q(y)(v)&=&(\mbox{Id}-\nabla\Pi(y))(v),\ \forall v \in\R^l\\
\nabla^2 Q(y)(v,w)&=&-\nabla^2\Pi (y)(v,w),\ \forall v,w \in\R^l.
\end{eqnarray*}
Set $A(y)(v,w)=-\nabla^2\Pi(y)(v,w)$ for $y\in N_{\delta_N}$ and $v,w\in\mathbb R^l$.
Then $F(u)$ can be rewritten by (see Gastel \cite{Gastel}):
$$ F(u)
=(-1)^{m+1}\sum\ls_{j=0}^{2m-2}\sum\ls_{|\alpha|=2m-2-j} \left(\begin{matrix}2m-2-j\\
\alpha\end{matrix}\right)\mbox{tr}^m(\nabla^j A)\circ
u(\nabla^{\alpha_1+1}u,\nabla^{\alpha_2+1}u,\cdots,\nabla^{\alpha_{j+2}+1}u).
$$
Direct calculations imply

\bex\begin{split}&\D^mQ(u)
=\nabla Q(u)(\D^m u)\\
+&\sum\ls_{j=0}^{2m-2}\sum\ls_{|\alpha|=2m-2-j} \left(\begin{matrix}2m-2-j\\
\alpha\end{matrix}\right)\mbox{trace}^m(\nabla^j A)\circ
u(\nabla^{\alpha_1+1}u,\nabla^{\alpha_2+1}u,\cdots,\nabla^{\alpha_{j+2}+1}u)\\
=& \nabla Q(u)(\D^m u)+(-1)^{m+1}F(u).
\end{split}\eex
Therefore we have

\beq\label{poly5.7} \begin{split} (\partial_t+(-1)^m\D^m)Q(u)
=&[\nabla Q(u)F(u)-F(u)]\\
=&-\nabla \Pi(u)(F(u)).
\end{split}\eeq
Multiplying both sides of (\ref{poly5.7}) by $Q(u)$ and integrating
over $\R^n$, we obtain that for $0\le t\le (\frac{R}{K_0})^{2m}$,

\beq\label{poly5.8}
\frac{d}{dt}\int_{\R^n}\rho(u)+\int_{\R^n}|\nabla^m Q(u)|^2
=-\int_{\R^n}\langle \nabla \Pi_(u)(F(u)), Q(u)\rangle=0,\eeq where
we have used the fact that $Q(u)\perp T_{{\Pi(u)}}N$ and $\nabla
\Pi(u)(F(u))\in T_{{\Pi(u)}}N$  on $\mathbb R^n\times [0,(\frac{R}{K_0})^{2m}]$
in the last step.

Since $\rho(u)|_{t=0}=0$, integrating (\ref{poly5.8}) with respect to $t$ implies 
$\rho(u)\equiv 0$ on $\R^n\times [0,(\frac{R}{K_0})^{2m}]$. Thus
$u(\R^n\times [0,(\frac{R}{K_0})^{2m}])\subset N$. Repeating the same argument also implies
that $u(\mathbb R^n\times [(\frac{R}{K_0})^{2m}, R^{2m}])\subset N$. 
This completes the proof of Theorem
\ref{main-theorem}. \endpf


\begin{thebibliography}{99}

\bibitem{Angelsberg-Pumberger}
G. Angelsberg, D. Pumberger, {\em A regularity result for polyharmonic maps with higher integrability}.
Ann. Global Anal. Geom.  35  (2009),  no. 1, 63-81.


\bibitem{Chen-Ding} Y. Chen, W. Ding, {\em Blow-up and global existence for heat
flows of harmonic maps.} Invent. Math. 99 (1990), no. 3, 567-578.

\bibitem{Chang-Ding-Ye} K. Chang, W. Ding, R. Ye, {\em Finite-time blow-up
of the heat flow of harmonic maps from surfaces.} J. Differen. Geom. 36(2),
507-515 (1992).

\bibitem{Coron-Ghidaglia} J. Coron, J. Ghidaglia, {\em Explosion en temps fini pour le
flot des applications harmoniques.} C.R. Acad. Sci. Paris 308, Serie I 339-344 (1989).

\bibitem{Chen-Lin} Y. Chen, F. Lin, {\em Evolution of harmonic maps with Dirichlet
boundary conditions.} Comm. Anal. Geom. 1 (1993), no. 3-4, 327-346..

\bibitem{Chen-Struwe} Y. Chen, M. Struwe, {\em Existence and partial regularity for
heat flow for harmonic maps.} Math. Z. 201, 83-103 (1989).

\bibitem{Eells-Sampson} J. Eells, J. Sampson, {\em Harmonic mappings of Riemannian
manifolds.} Amer. J. Math. 86, 109-160 (1964)

\bibitem{Gastel} A. Gastel, {\em The extrinsic polyharmonic map heat flow in the
critical dimension.} Adv. Geom. 6 (2006), no. 4, 501-521.

\bibitem{Gastel-Scheven}
A. Gastel, C. Scheven, {\em Regularity of polyharmonic maps in the critical dimension.}
Comm. Anal. Geom.  17  (2009),  no. 2, 185-226.

\bibitem{GSZ} P. Goldstein, P. Strzelecki, A. Zatorska-Goldstein, 
{\em On polyharmonic maps into spheres in the critical dimension}.  
Ann. Inst. H. Poincar\'e Anal. Non Lin\'eaire,
vol. 26, issue 4, pp. 1387-1405.

\bibitem{Hildebrandt-Kaul-Widman} S. Hildebrandt, H. Kaul, K. Widman,
{\em An existence theorem for harmonic mappings of Riemannian manifolds.} Acta Math.
138 (1977), no. 1-2, 1-16.

\bibitem{Hormander} L. H$\ddot{\mbox{o}}$rmander,
The analysis of linear partial differential operators I, 2nd edition,
Springer-Verlag Berlin Heidelberg, 1990.

\bibitem{Koch-Lamm} H. Koch, T. Lamm, {\em Geometric flows with rough initial data.}
arXiv: 0902.1488v1, 2009.

\bibitem{Koch-Tataru} H. Koch, D. Tataru, {\em Well-posedness for the Navier-Stokes
equations.} Adv. Math. 157 (2001), no. 1, 22-35.

\bibitem{Lamm2001} T. Lamm, Biharmonischer W$\ddot{\mbox{a}}$rmefluss.
Diplomarbeit Universit$\ddot{\mbox{a}}$t Freiburg (2001).

\bibitem{Lamm2004} T. Lamm, {\em Heat flow for extrinsic biharmonic maps
with small initial energy.} Ann. Global Anal. Geom. 26 (2004), no. 4, 369-384.

\bibitem{Lamm2005} T. Lamm, {\em Biharmonic map heat flow into manifolds of
nonpositive curvature.} Calc. Var. Partial Differential Equations 22 (2005),
no. 4, 421-445.

\bibitem{Lamm-Wang} T. Lamm, C. Wang, {\em Boundary regularity for polyharmonic maps in the critical dimension.}
Adv. Calc. Var.  2  (2009),  no. 1, 1-16.


\bibitem{Lin-Wang} F. Lin, C. Wang, The analysis of harmonic maps and their heat
flows. World Scientific Publishing Co. Pte. Ltd., Hackensack, NJ, 2008.


\bibitem{Moser} R. Moser, {\em Regularity of minimizing extrinsic polyharmonic maps in the critical dimension}.
Preprint (2009).

\bibitem{Stein} E. Stein, Harmonic analysis, Vol. 43 of Princeton
Mathematical Series, Princeton University Press, Princeton, NJ,
1993.

\bibitem{Struwe} M. Struwe, {\em On the evolution of harmonic maps of
Riemannian surfaces.} Comment. Math. Helv. 60, 558-581 (1985).

\bibitem{Wang2007} C. Wang, {\em Heat flow of biharmonic maps in dimensions
four and its application}. Pure Appl. Math. Q. 3 (2007), no. 2, part 1, 595-613.



\bibitem{Wang-biharmonic} C. Wang,{\em Well-posedness for the heat flow
of biharmonic maps with rough initial data.} Preprint (2010).

\bibitem{Wang-harmonic} C. Wang, {\em Well-posedness for the heat flow of
harmonic maps and the liquid crystal flow with rough initial data.} Preprint (2010).

\end{thebibliography}
\end{document}